\documentclass{amsart}
\usepackage{latexsym,amsmath,amssymb}
\theoremstyle{definition}
\theoremstyle{remark}

\numberwithin{equation}{section}
\begin{document}

\title[Weighted conditional type operators ]
{ Paranormal weighted conditional type operators }

\author{\sc\bf Y. Estaremi  }
\address{\sc Y. Estaremi  }
\email{yestaremi@pnu.ac.ir}

\address{Department of mathematics, Payame Noor university , p. o. box: 19395-3697, Tehran,
Iran}

\thanks{}

\thanks{}

\subjclass[2010]{47B20}

\keywords{Conditional expectation, paranormal operators, quasi-$\ast$-paranormal operators,
spectrum, point
spectrum, approximate point spectrum. }

\date{}

\dedicatory{}

\commby{}

%%% ----------------------------------------------------------------------
\begin{abstract}
In this paper, some sub-classes of paranormal weighted conditional
expectation type operators, such as $\ast$-paranormal, quasi-$\ast$-paranormal
and $(n,k)$-quasi-$\ast$-paranormal weighted conditional
expectation type operators on $L^{2}(\Sigma)$ are investigated.
Also, some applications about the spectrum, point spectrum, joint point spectrum, approximate point spectrum and joint approximate point spectrum of these classes are presented.

\noindent {}
\end{abstract}

\maketitle

\section{ \sc\bf Introduction and Preliminaries}

Theory of weighted conditional
expectation type operators is one of important
arguments in the connection of operator theory and measure theory.
Weighted conditional expectations have been studied in an operator
theoretic setting, by many authors, for example, De pagter and Grobler \cite{g}
and Rao \cite{rao1, rao2}, as positive operators acting on
$L^p$-spaces or Banach function spaces. In \cite{mo}, S.-T. C. Moy
characterized all operators on $L^p$ of the form $f\rightarrow
E(fg)$ for $g$ in $L^q$ with $E(|g|)$ bounded.  Also, some results
about these operators can be
found in \cite{dou, her, lam}. In \cite{dhd} P.G. Dodds, C.B.
Huijsmans and B. de Pagter showed that lots of operators are of
the form of weighted conditional type operators. In \cite{g} a class of operators which factorizes
through weighted conditional type operators is investigated. This class of operators
includes operators such as kernel operators and order continuous
Riesz homomorphisms. Also, we investigated some classical properties of these operators on $L^p$-spaces in \cite{e1,e2,ej}.\\

\vspace*{0.3cm} Let $(X,\Sigma,\mu)$ be a $\sigma$-finite measure space. For a sub-$\sigma$-finite algebra
$\mathcal{A}\subseteq\Sigma$, the conditional expectation operator
associated with $\mathcal{A}$ is the mapping $f\rightarrow
E^{\mathcal{A}}f$, defined for all non-negative measurable
functions $f$ as well as for all $f\in L^2(\Sigma)$, where
$E^{\mathcal{A}}f$, by the Radon-Nikodym theorem, is the unique
$\mathcal{A}$-measurable function satisfying
$\int_{A}fd\mu=\int_{A}E^{\mathcal{A}}fd\mu, \ \ \ \forall A\in
\mathcal{A}$. As an operator on $L^{2}({\Sigma})$,
$E^{\mathcal{A}}$ is idempotent and
$E^{\mathcal{A}}(L^2(\Sigma))=L^2(\mathcal{A})$. This operator
will play a major role in our work. Let $f\in L^0(\Sigma)$, then
$f$ is said to be conditionable with respect to $E$ if
$f\in\mathcal{D}(E):=\{g\in L^0(\Sigma): E(|g|)\in
L^0(\mathcal{A})\}$ in which $L^0(\Sigma)$ is the vector space of all equivalence classes of
almost everywhere finite valued measurable functions on $X$. Throughout this paper we take $u$ and $w$ in
$\mathcal{D}(E)$. If there is no possibility of confusion, we
write $E(f)$ in place of $E^{\mathcal{A}}(f)$.  A detailed
discussion about this operator may be found in \cite{rao}. All
comparisons between two functions or two sets are to be
interpreted as holding up to a $\mu$-null set. The support of a
measurable function $f$ is defined as $S(f)=\{x\in X; f(x)\neq
0\}$.

Let $\mathcal{H}$ be the infinite dimensional complex Hilbert
space and let $\mathcal{L(H)}$ be the algebra of all bounded
operators on $\mathcal{H}$. An operator $T\in \mathcal{L(H)}$ is a
partial isometry if $\|Th\|=\|h\|$ for $h$ orthogonal to the
kernel of $T$. It is known that an operator $T$ on a Hilbert space
is partial isometry if and only if $TT^{\ast}T=T$. Every operator
$T$ on a Hilbert space $\mathcal{H}$ can be decomposed into $T =
U|T|$ with a partial isometry $U$, where $|T| =
(T^*T)^{\frac{1}{2}}$ . $U$ is determined uniquely by the kernel
condition $\mathcal{N}(U) = \mathcal{N}(|T|)$. Then this
decomposition is called the polar decomposition. The Aluthge
transformation $\widehat{T}$ of the operator $T$ is defined by
$\widehat{T}=|T|^{\frac{1}{2}}U|T|^{\frac{1}{2}}$. The operator
$T$ is said to be a positive operator and written as $T\geq 0$, if
$\langle
Th, h\rangle\geq 0$, for all $h\in \mathcal{H}$. \\

In this paper we will be concerned with characterizing weighted
conditional expectation type operators on $L^2(\Sigma)$ in terms
of membership of the partial paranormal classes. Also, we prove that the
point spectrum and joint point spectrum of the weighted conditional type operators $M_wEM_u$ are the same, when $u,w$ satisfy a mild condition. Here is
a brief review of what constitutes membership for an operator $T$
on a Hilbert space in each classes:\\

(i) $T$ is called paranormal if for all unit vectors $x$ in
$H$, $\|Tx\|^2\leq \|T^2x\|$ or equivalently
$$T^{\ast^2}T^2-2\lambda T^{\ast}T+\lambda^2I\geq0,$$
for all $\lambda>0$;\\

(ii) $T$ is called $M$-paranormal, if  there exists $M>0$ such that for all unit vectors $x$ in $H$,
$\|Tx\|^2\leq M\|T^2x\|$ or equivalently $M^2T^{\ast^2}T^2-2\lambda T^{\ast}T+\lambda^2I\geq0$, for all $\lambda>0$.\\

(iii) $T$ is called $\ast$-paranormal, if $\|T^{\ast}x\|^2\leq \|T^2x\|$ for all unit vector $x\in H$ or equivalently
$$T^{\ast^2}T^2-2\lambda TT^{\ast}
+\lambda^2I\geq0,$$
for all $\lambda>0$;\\

(iv) $T$ is called quasi-$\ast$-paranormal, if it satisfies the following inequality:
$$\|T^{\ast}Tx\|^2\leq\|T^3x\|.\|Tx\|$$
for all $x\in H$  or equivalently
$$T^{\ast}(T^{\ast^2}T^2-2\lambda TT^{\ast}
+\lambda^2I)T\geq0,$$
for all $\lambda>0$;\\

(v) $T$  is called $(n,k)$-quasi-$\ast$-paranormal if
$$\|T^{1+n}(T^kx)\|^{\frac{1}{1+n}}\|T^kx\|^{\frac{n}{n+1}}\geq\|T^{\ast}T^kx)\|$$
for all $x\in H$, or equivalently
$$T^{\ast^k}T^{\ast^{(1+n)}}T^{1+n}T^k-(1+n)\mu^nT^{\ast^k}TT^{\ast}T^k+n\mu^{n+1}T^{\ast^k}T^k\geq0$$
for all $\mu>0$;\\

(vi) $T$ is called absolute-$k$-paranormal for each $k>0$ if
$$\||T|^kTx\|\geq \|Tx\|^{k+1}$$
for every unit vector $x\in H$ or equivalently
$$T^{\ast}|T|^{2k}T-(k+1)\lambda^k|T|^2+k\lambda^{k+1}I\geq0$$
for all $\lambda>0$;\\

\section{ \sc\bf Partial paranormal weighted conditional type operators}
First we recall some results of \cite{dhd} that state our results is valid for a large class of linear operators. Let $(X, \Sigma, \mu)$ be a finite measure space, then $L^{\infty}(\Sigma)\subseteq L^2(\Sigma)\subseteq L^1(\Sigma)$ and $L^2(\Sigma)$ is an order ideal of measurable functions on $(X,\Sigma,\mu)$. Thus by propositions $(3.1, 3.3, 3.6)$ of \cite{dhd} we have theorems A, B, C:\\

\vspace*{0.4cm} {\bf Theorem A.} If $T$ is a linear operator on $L^2(\Sigma)$ for which\\

(i) $Tf\in L^{\infty}(\Sigma)$ whenever $f\in L^{\infty}(\Sigma)$.\\

(ii) $\|Tf_n\|_1\rightarrow 0$ for all sequences $\{f_n\}_{n=1}^{\infty}\subseteq L^{2}(\Sigma)$ such that $|f_n|\leq g$ $(n=1,2,3,....)$ for some $g\in L^{2}(\Sigma)$ and $f_n\rightarrow 0$ a.e.,\\

(iii) $T(f.Tg)=Tf.Tg$ for all $f\in L^{\infty}(\Sigma)$ and all $g\in L^{2}(\Sigma)$,\\

then there exists a $\sigma$-subalgebra $\mathcal{A}$ of $\Sigma$ and there exists $w\in L^{2}(\Sigma)$ such that $Tf=E^{\mathcal{A}}(wf)$ for all $f\in L^{2}(\Sigma)$.\\

 \vspace*{0.4cm} {\bf Theorem B.} For a linear operator $T:L^2(\Sigma)\rightarrow L^2(\Sigma)$ the following statement are equivalent.\\

(i) $T$ is positive and order continuous, $T^2=T$, $T1=1$ and the range $\mathcal{R}(T)$ of $T$ is a sublattice.\\

(ii) There exist a $\sigma$-subalgebra $\mathcal{A}$ of $\Sigma$ and a function $0\leq w\in L^{2}(\Sigma)$ with $E^{\mathcal{A}}(w)=1$  such that $Tf=E^{\mathcal{A}}(wf)$ for all $f\in L^{2}(\Sigma)$.\\

\vspace*{0.4cm} {\bf Theorem C.} For a linear operator $T:L^2(\Sigma)\rightarrow L^2(\Sigma)$ the following statement are equivalent.\\

(i) $T$ is a positive and order continuous projection onto a sublattice such that $T1$ is strictly positive.\\

(ii) There exist a $\sigma$-subalgebra $\mathcal{A}$ of $\Sigma$, $0\leq w\in L^{2}(\Sigma)$  and a strictly positive function $k\in L^1(\Sigma)$ with $E^{\mathcal{A}}(wk)=1$  such that $Tf=E^{\mathcal{A}}(wf)$ for all $f\in L^{2}(\Sigma)$. Moreover, if we choose $k$ such that $E^{\mathcal{A}}(k)=1$, then both $w$ and $k$ are uniquely determined by $T$. \\

Here, we recall some properties of weighted conditional
type operators, that we have proved in
\cite{ej}.\\

 The operator $T=M_wEM_u$ is bounded on $L^{2}(\Sigma)$ if and
only if $(E|w|^{2})^{\frac{1}{2}}(E|u|^{2})^{\frac{1}{2}} \in
L^{\infty}(\mathcal{A})$, and in this case its norm is given by
$\|T\|=\|(E(|w|^{2}))^{\frac{1}{2}}(E(|u|^{2}))^{\frac{1}{2}}\|_{\infty}$.
The unique polar decomposition of bounded operator $T=M_wEM_u$ is
$U|T|$, where

$$|T|(f)=\left(\frac{E(|w|^{2})}{E(|u|^{2})}\right)^{\frac{1}{2}}\chi_{S}\bar{u}E(uf)$$
and
 $$U(f)=\left(\frac{\chi_{S\cap
 G}}{E(|w|^{2})E(|u|^{2})}\right)^{\frac{1}{2}}wE(uf),$$
for all $f\in L^{2}(\Sigma)$, where $S=S(E(|u|^2))$ and
$G=S(E(|w|^2))$.\\

In the sequel some necessary and sufficient conditions for
weighted conditional type operator $M_wEM_u$ to be $M$-paranormal, quasi-$\ast$-paranormal, absolute-$k$-paranormal and $(n,k)$-quasi-$\ast$-paranormal will be presented.\\

\vspace*{0.4cm} {\bf Theorem 2.1.} Let $T=M_wEM_u$ be a bounded operator on $L^2(\Sigma)$. Then\\

(i) If $T$ is $M$- paranormal, we have
$$(M^2|E(uw)|^2E(|w|^2)-2\lambda E(|w|^2))|E(u)|^2+\lambda^2\geq0,$$
for all $\lambda>0$.\\

(ii) If $M^2|E(uw)|^2E(|w|^2)-2\lambda E(|w|^2)\geq0$ for all $\lambda>0$, then $T$ is $M$-paranormal.\\

\vspace*{0.3cm} {\bf Proof.} (i)
By induction and by Lemma we have
\begin{align*}
T^{\ast}T&=M_{E(|w|^2)}M_{\bar{u}}EM_{u},\\
T^{\ast^2}T^2&=M_{|E(uw)|^{2}E(|w|^2)}M_{\bar{u}}EM_{u}.
\end{align*}
So for every $\lambda>0$ and $M>0$
\begin{align*}
M^2T^{\ast^2}T^2-2\lambda T^{\ast}T+\lambda^2I&=M^2M_{|E(uw)|^{2}E(|w|^2)}M_{\bar{u}}EM_{u}-2\lambda M_{E(|w|^2)}M_{\bar{u}}EM_{u}+\lambda^2I\\
&=(M_{M^2|E(uw)|^{2}E(|w|^2)-2\lambda E(|w|^2)})M_{\bar{u}}EM_{u}+\lambda^2I.
\end{align*}
Let $\alpha=M^2|E(uw)|^{2}E(|w|^2)-2\lambda E(|w|^2)$. Then for every $f\in L^2(\Sigma)$ we get

\begin{align*}
\langle M_{\alpha}M_{\bar{u}}EM_{u}f+\lambda^2f,f\rangle&=\int_{X}\alpha|E(uf)|^2d\mu+\int_{X}\lambda^2|f|^2d\mu
\end{align*}

This implies that if $\alpha\geq0$ a.e, then $T$ is $M$-paranormal.\\

(ii) If $T$ is $M$-paranormal, then for all $f\in L^2(\mathcal{A})$

\begin{align*}
\langle M_{\alpha}M_{\bar{u}}EM_{u}f+\lambda^2f,f\rangle&=\int_{X}\alpha|E(uf)|^2d\mu+\int_{X}\lambda^2|f|^2d\mu\\
&=\int_{X}\alpha|E(u)|^2|f|^2d\mu+\int_{X}\lambda^2|f|^2d\mu\\
&=\int_{X}(\alpha|E(u)|^2+\lambda^2)|f|^2d\mu\geq0.\\
\end{align*}
Therefore $\alpha|E(u)|^2+\lambda^2\geq0$ a.e,.\\

\vspace*{0.4cm} {\bf Corollary 2.2.} Let $T=M_wEM_u$ be a bounded operator on $L^2(\Sigma)$. Then\\

(i) If $T$ is paranormal, we have
$$(|E(uw)|^2E(|w|^2)-2kE(|w|^2))|E(u)|^2+k^2\geq0,$$

(ii) If $|E(uw)|^2E(|w|^2)-2kE(|w|^2)\geq0$, then $T$ is paranormal.\\

\vspace*{0.4cm} {\bf Corollary 2.3.} Let $u\in L^0(\mathcal{A})$ and $T=M_wEM_u$ be a bounded operator on $L^2(\Sigma)$. Then:\\

 (i) $T$ is $M$-paranormal if and only if $$(M^2|uE(w)|^2E(|w|^2)-2\lambda E(|w|^2))|u|^2+\lambda^2\geq0,$$
for all $\lambda>0$.\\

(ii) $T$ is paranormal if and only if $$(|uE(w)|^2E(|w|^2)-2\lambda E(|w|^2))|u|^2+\lambda^2\geq0,$$
for all $\lambda>0$.\\

\vspace*{0.3cm} {\bf Proof.} Since $|E(F)|^2\leq E(|f|^2)$ for every $f\in L^2(\Sigma)$, then by similar method of Theorem 2.2 we get the proof.\\

The definition of quasi-$\ast$-paranormal and $\ast$-paranormal operators shows that, if $T$ is quasi-$\ast$-paranormal, then $T\mid_{\overline{\mathcal{R}(T)}}$ is $\ast$-paranormal. Therefore, if $T$ has dense range, then $T$ is quasi-$\ast$-paranormal if and only if is $\ast$-paranormal. In the next theorem we give a necessary and sufficient condition for $M_wEM_u$ to be quasi-$\ast$-paranormal.\\

\vspace*{0.4cm} {\bf Theorem 2.4.} Let $T=M_wEM_u$ be a bounded operator on $L^2(\Sigma)$. Then $T$ is quasi-$\ast$-paranormal if and only if
$$E(|u|^2)E(|w|^{2})\leq|E(uw)|^{2} a.e,.\ \ \ on \ \ G.$$
Where $G=S(E(|w|^2)$.\\

\vspace*{0.3cm} {\bf Proof.} By direct computations we have
$$T^nf=(E(uw))^{n-1}wE(uf),  \ \ \ \ T^{\ast ^n}f=(\overline{E(uw)})^{n-1}\bar{u}E(\bar{w}f),$$

 for all $f\in L^2(\Sigma)$ and $n\in \mathbb{N}$. So we get that
 \begin{align*}
 (T^{\ast}T)^2&=M_{\bar{u}E(|u|^{2})\chi_{S}(E(|w|^{2}))^{2}}EM_{u},\\
T^{\ast}T&=M_{\bar{u}E(|w|^{2})}EM_{u},\\
T^{\ast^3}T^3&=M_{|E(uw)|^{4}E(|w|^2)}M_{\bar{u}}EM_{u}.
 \end{align*}
 Therefore $T$ is quasi-$\ast$-paranormal if and only if
$$
 M_{|E(uw)|^{4}E(|w|^2)}M_{\bar{u}}EM_{u}-2\lambda M_{\bar{u}E(|u|^{2})\chi_{S}(E(|w|^{2}))^{2}}EM_{u}
 +\lambda^2 M_{\bar{u}E(|w|^{2})}EM_{u}$$
 $$=\left(M_{|E(uw)|^{4}E(|w|^2)}
 -2\lambda M_{E(|u|^{2})\chi_{S}(E(|w|^{2}))^{2}}+\lambda^2 M_{E(|w|^{2})}\right)M_{\bar{u}}EM_{u} \geq0.$$

This implies that $T$ is quasi-$\ast$-paranormal if and only if for all $f\in L^2(\Sigma)$ and $\lambda>0$
\begin{align*}
0\leq&\langle M_{|E(uw)|^{4}E(|w|^2)}f -2\lambda M_{E(|u|^{2})\chi_{S}(E(|w|^{2}))^{2}}f+\lambda^2 M_{E(|w|^{2})}f,f\rangle\\
&=\langle M_{\left(|E(uw)|^{4}E(|w|^2) -2\lambda E(|u|^{2})\chi_{S}(E(|w|^{2}))^{2}+\lambda^2 E(|w|^{2})\right)}f,f\rangle,\\
&\Longleftrightarrow M_{\left(|E(uw)|^{4}E(|w|^2) -2\lambda E(|u|^{2})\chi_{S}(E(|w|^{2}))^{2}+\lambda^2 E(|w|^{2})\right)}\geq0\\
&\Longleftrightarrow|E(uw)|^{4}E(|w|^2) -2\lambda E(|u|^{2})\chi_{S}(E(|w|^{2}))^{2}+\lambda^2 E(|w|^{2})\geq0\\
&\Longleftrightarrow (E(|u|^2))^2(E(|w|^{2}))^{4}-(E(|w|^{2}))^2|E(uw)|^{4}\geq0\\
&\Longleftrightarrow E(|u|^2)E(|w|^{2})\leq|E(uw)|^{2} \ \  on \ G,
\end{align*}

where we have used the fact that $T_{1}T_{2}\geq0$ if $T_{1}\geq0$,
$T_{2}\geq0$ and $T_{1}T_{2}=T_{2}T_{1}$ for all $T_{i}\in
\mathcal{B}(\mathcal{H})$, positivity of $M_{\bar{u}}EM_{u}$ and $M_{\alpha}M_{\bar{u}}EM_{u}=M_{\bar{u}}EM_{u}M_{\alpha}$ such that
$$\alpha=|E(uw)|^{4}E(|w|^2) -2\lambda E(|u|^{2})\chi_{S}(E(|w|^{2}))^{2}+\lambda^2 E(|w|^{2}).$$

An operator $T\in \mathcal{L(H)}$ is a quasi-$\ast$-$A$-class operator if
$T^{\ast}|T^2|T\geq T^{\ast}|T^{\ast}|^2T$. And $T$ is an $A$-class operator if $|T|^2\leq|T^2|$. One can see \cite{fu} for more details. In \cite{e2} we studied quasi-$\ast$-$A$-class weighted conditional type operators. Here we get that quasi-$\ast$-$A$-class and quasi-$\ast$-paranormal weighted conditional type operators are coincided.\\

\vspace*{0.4cm} {\bf Theorem 2.5.} Let $T=M_wEM_u$ be a bounded operator on $L^2(\Sigma)$ and $G=X$. Then the followings are mutually equivalent:\\

(a) $T$ is quasi-$\ast$-paranormal;\\

(b) $T$ is a quasi-$\ast$-$A$-class operator;\\

(c) $E(|u|^2)E(|w|^{2})\leq|E(uw)|^{2}$ a.e,.\\

Moreover, if $S(E(u))=X=G$, then (a), (b), (c) and (d) are mutually equivalent, where\\

(d) $T$ is an $A$-class operator.\\

\vspace*{0.4cm} {\bf Proof.} This is a direct consequence of Theorem 2.4 and Theorems 2.6 and 2.8 of  \cite{e2}.\\

\vspace*{0.4cm} {\bf Theorem 2.6.} Let $T=M_wEM_u$ be a bounded operator on $L^2(\Sigma)$. Then,\\

(i) If $T$ is absolute-$k$-paranormal, we have
$$(|E(uw)|^2(E(|u|^2))^{k-1}.\chi_S(E(|w|^2))^k|E(u)|^2-(k+1)\lambda^kE(|w|^2))|E(u)|^2+k\lambda^{k+1}\geq0,$$
for all $\lambda>0$.\\

(ii) If $|E(uw)|^2(E(|u|^2))^{k-1}.\chi_S(E(|w|^2))^k-(k+1)\lambda^kE(|w|^2)\geq0$ for all $\lambda>0$, then $T$ is absolute-$k$-paranormal.\\

\vspace*{0.3cm} {\bf Proof.} (i) By similar methods of last theorems we have
\begin{align*}
T^{\ast}|T|^{2k}T&=M_{(E(|u|^{2}))^{k-1}\chi_{S}(E(|w|^{2}))^{k}|E(uw)|^2}M_{\bar{u}}EM_{u},\\
T^{\ast}T&=M_{E(|w|^{2})}M_{\bar{u}}EM_{u},\\
\end{align*}
If $T$ is absolute-$k$-paranormal, then for all $f\in L^2(\mathcal{A})$
\begin{align*}
0&\leq\langle M_{(E(|u|^{2}))^{k-1}\chi_{S}(E(|w|^{2}))^{k}|E(uw)|^2}M_{\bar{u}}EM_{u}f-(k+1)\lambda^kM_{E(|w|^{2})}M_{\bar{u}}EM_{u}f+k\lambda^{k+1}f,f\rangle\\
&=\int_{X}\left((E(|u|^{2}))^{k-1}\chi_{S}(E(|w|^{2}))^{k}|E(uw)|^2-(k+1)\lambda^kE(|w|^{2})\right)\bar{u}E(u)f\bar{f}d\mu+k\lambda^{k+1}\int_{X}|f|^2d\mu\\
&=\int_{X}\left((E(|u|^{2}))^{k-1}\chi_{S}(E(|w|^{2}))^{k}|E(uw)|^2|E(u)|^2-(k+1)\lambda^kE(|w|^{2})|E(u)|^2
+k\lambda^{k+1}\right)|f|^2d\mu.
\end{align*}
So we get that
$$(E(|u|^{2}))^{k-1}\chi_{S}(E(|w|^{2}))^{k}|E(uw)|^2|E(u)|^2-(k+1)\lambda^kE(|w|^{2})|E(u)|^2
+k\lambda^{k+1}\geq0.$$

(ii) The operator $T$ is absolute-$k$-paranormal if for all $f\in L^2(\Sigma)$

\begin{align*}
0&\leq T^{\ast}|T|^{2k}T-(k+1)\lambda^k|T|^2+k\lambda^{k+1}I\geq0\\
&=\langle M_{(E(|u|^{2}))^{k-1}\chi_{S}(E(|w|^{2}))^{k}|E(uw)|^2}M_{\bar{u}}EM_{u}f-(k+1)\lambda^kM_{E(|w|^{2})}M_{\bar{u}}EM_{u}f+k\lambda^{k+1}f,f\rangle\\
&=\int_{X}\left((E(|u|^{2}))^{k-1}\chi_{S}(E(|w|^{2}))^{k}|E(uw)|^2-(k+1)\lambda^kE(|w|^{2})\right)|E(uf)|^2d\mu+k\lambda^{k+1}\int_{X}|f|^2d\mu.\\
\end{align*}
This implies that if $|E(uw)|^2(E(|u|^2))^{k-1}.\chi_S(E(|w|^2))^k-(k+1)\lambda^kE(|w|^2)\geq0$, then $T$ is absolute-$k$-paranormal.\\

\vspace*{0.4cm} {\bf Corollary 2.7.} Let $u\in L^0(\mathcal{A})$ and $T=M_wEM_u$ be a bounded operator on $L^2(\Sigma)$. Then $T$ is absolute-$k$-paranormal if and only if $$(|uE(w)|^2|u|^{2k-2}.\chi_S(E(|w|^2))^k|u|^2-(k+1)\lambda^kE(|w|^2))|u|^2+k\lambda^{k+1}\geq0,$$
for all $\lambda>0$.\\

\vspace*{0.3cm} {\bf Proof.} Since $|E(F)|^2\leq E(|f|^2)$ for every $f\in L^2(\Sigma)$, then by similar method of Theorem 2.6 we get the proof.\\

\vspace*{0.4cm} {\bf Theorem 2.8.} Let $T=M_wEM_u$ be a bounded operator on $L^2(\Sigma)$. Then $T$ is $(n,k)$-quasi-$\ast$-paranormal if and only if $\alpha_{\mu}\geq0$ for all $\mu>0$, where
\begin{align*}
\alpha_k=|E(uw)|^{2(n+k)}E(|w|^2)&-(1+n)\mu^n|E(uw)|^{2(k-1)}\chi_{S_0}(E(|w|^2))^2E(|u|^2)\\
&+\mu^{1+n}|E(uw)|^{2(k-1)}\chi_{S_0}E(|w|^2),
\end{align*}
and $S_0=S(E(uw))$.\\

\vspace*{0.3cm} {\bf Proof.} By similar methods of last theorems we have

 \begin{align*}
 T^{\ast^k}T^{\ast^{(1+n)}}T^{1+n}T^k&=M_{|E(uw)|^{2(n+k)}E(|w|^2)}M_{\bar{u}}EM_{u},\\
T^{\ast^k}TT^{\ast}T^k&=M_{|E(uw)|^{2(k-1)}(E(|w|^2))^2E(|u|^2)}M_{\bar{u}}EM_{u},\\
T^{\ast^k}T^k&=M_{|E(uw)|^{2(k-1)}E(|w|^2)}M_{\bar{u}}EM_{u}.
 \end{align*}
 This implies that $T$ is $(n,k)$-quasi-$\ast$-paranormal if and only if
\begin{align*}
 (M_{|E(uw)|^{2(n+k)}E(|w|^2)}&-(1+n)\mu^nM_{|E(uw)|^{2(k-1)}(E(|w|^2))^2E(|u|^2)}\\
 &+n\mu^{n+1}M_{|E(uw)|^{2(k-1)}E(|w|^2)})M_{\bar{u}}EM_{u}\geq0.
\end{align*}
This inequality holds if and only if

\begin{align*}
 M_{|E(uw)|^{2(n+k)}E(|w|^2)}&-(1+n)\mu^nM_{|E(uw)|^{2(k-1)}(E(|w|^2))^2E(|u|^2)}\\
 &+n\mu^{n+1}M_{|E(uw)|^{2(k-1)}E(|w|^2)}\geq0.
\end{align*}

where we have used the fact that $T_{1}T_{2}\geq0$ if $T_{1}\geq0$,
$T_{2}\geq0$ and $T_{1}T_{2}=T_{2}T_{1}$ for all $T_{i}\in
\mathcal{B}(\mathcal{H})$, positivity of $M_{\bar{u}}EM_{u}$ and $M_{\alpha_{\mu}}M_{\bar{u}}EM_{u}=M_{\bar{u}}EM_{u}M_{\alpha_{\mu}}$ such that
\begin{align*}
\alpha_{\mu}=|E(uw)|^{2(n+k)}E(|w|^2)&-(1+n)\mu^n|E(uw)|^{2(k-1)}\chi_{S_0}(E(|w|^2))^2E(|u|^2)\\
&+n\mu^{n+1}|E(uw)|^{2(k-1)}\chi_{S_0}E(|w|^2).
\end{align*}

 Therefore the operator $T$ is $(n,k)$-quasi-$\ast$-paranormal if and only if $M_{\alpha_{\mu}}\geq0$ if and only if $\alpha_{\mu}\geq0$ for all $\mu>0$.\\

Recall that an operator $T\in \mathcal{L(H)}$ is said to be $n$-$\ast$-paranormal if $\|T^{n+1}x\|^{\frac{1}{n+1}}\|x\|^{\frac{n}{n+1}}\geq \|T^{\ast}x\|$ for all $x\in \mathcal{H}$. Also, $T$ is called $k$-quasi -$\ast$-paranormal if $\|T^{2}(T^kx)\|^{\frac{1}{2}}\|T^kx\|^{\frac{1}{2}}\geq \|T^{\ast}(Tx)\|$ for all $x\in \mathcal{H}$. So, we get the following corollaries.\\

\vspace*{0.4cm} {\bf Corollary 2.9.} Let $T=M_wEM_u$ be a bounded operator on $L^2(\Sigma)$. Then $T$ is $n$-$\ast$-paranormal if and only if $\alpha_{\mu}\geq0$ for all $\mu>0$, where

\begin{align*}
\alpha_{\mu}=|E(uw)|^{2(n)}E(|w|^2)&-(1+n)\mu^n|E(uw)|^{-2}\chi_{S_0}(E(|w|^2))^2E(|u|^2)\\
&+n\mu^{n+1}|E(uw)|^{-2}\chi_{S_0}E(|w|^2).
\end{align*}

\vspace*{0.4cm} {\bf Corollary 2.10.} Let $T=M_wEM_u$ be a bounded operator on $L^2(\Sigma)$. Then $T$ is $k$-quasi-$\ast$-paranormal if and only if $\alpha_{\mu}\geq0$ for all $\mu>0$, where
\begin{align*}
\alpha_{\mu}=|E(uw)|^{2(1+k)}E(|w|^2)&-2\mu^n|E(uw)|^{2(k-1)}\chi_{S_0}(E(|w|^2))^2E(|u|^2)\\
&+\mu^{2}|E(uw)|^{2(k-1)}\chi_{S_0}E(|w|^2).
\end{align*}
%By Theorem 2.3 we have
%
%$$|T|^2(f)=E(|w|^2)\chi_{S}\bar{u}E(uf), \ \ \ \
%|T^2|(f)=|E(uw)|(\frac{E(|w|^2)}{E(|u|^2)})^{\frac{1}{2}}\chi_{S}\bar{u}E(uf),$$
%for all $f\in L^2(\Sigma)$. So for every $f\in L^2(\Sigma)$ we
%have
%
%$$\langle|T^2|(f)-|T|^2(f),f\rangle=\int_{X}|E(uw)|(\frac{E(|w|^2)}{E(|u|^2)})^{\frac{1}{2}}\chi_{S}\overline{uf}E(uf)-E(|w|^2)\chi_{S}\overline{uf}E(uf)d\mu$$
%$$=\int_{X}|E(uw)|(\frac{E(|w|^2)}{E(|u|^2)})^{\frac{1}{2}}\chi_{S}|E(uf)|^2-E(|w|^2)\chi_{S}|E(uf)|^2d\mu.$$

\section{ \sc\bf Some Applications}

For $T\in \mathcal{L(H)}$, let $\sigma_p(T)$, $\sigma_{jp}(T)$, $\sigma_{a}(T)$ and $\sigma_{ja}(T)$ denote the point spectrum, joint point spectrum, approximate point spectrum and joint approximate point spectrum of $T$. $T$ is called isoloid if every isolated point of $\sigma(T)$ is an eigenvalue of $T$. Let $\mu\in \mathbb{C}$ be an isolated point of $\sigma(T)$. Then the riesz idempotent $E_{\mu}$ of $T$ with respect to $\mu$ is defined by
$$E_{\mu}:=\frac{1}{2\pi i}\int_{\partial D_{\mu}}(\mu I-T)^{-1}d\mu,$$
where $D_{\mu}$ is the closed disk centered at $\mu$ which contains no other points of $\sigma(T)$.\\

An operator $T\in \mathcal{L(H)}$ is said to have single valued extension property at $\lambda_0\in \mathbb{C}$ (SVEP at $\lambda_0$ for brevity) if for every open neighborhood $U$ of $\lambda_0$, the only analytic function $f:U\rightarrow H$ which satisfies the equation $(\lambda I-T)f(\lambda)=0$ for all$\lambda\in U$ is the constant function $f\equiv 0$. Here we recall some spectral results about $M_wEM_u$. Then we get some conclusions for quasi-$\ast$-paranormal and \\

\vspace*{0.3cm} {\bf Theorem 3.1. \cite{e2}} Let
$T=M_wEM_u:L^2(\Sigma)\rightarrow L^2(\Sigma)$. Then\\

(a) $$\sigma(M_wEM_u)\setminus \{0\}=ess \
range(E(uw))\setminus\{0\}.$$

(b) If $S\cap G=X$, then $$\sigma(M_wEM_u)=ess \ range(E(uw)),$$
where $S=S(E(|u|^2))$ and $G=S(E(|w|^2))$.

(c)
$$\sigma_{p}(M_wEM_u)\setminus\{0\}=\{\lambda\in\mathbb{C}\setminus\{0\}:\mu(A_{\lambda,w})>0\},$$
where $A_{\lambda,w}=\{x\in X:E(uw)(x)=\lambda\}$.\\

\vspace*{0.3cm} {\bf Proposition 3.2,} If
$|E(uw)|^2\geq E(|u|^2)E(|w|^2)$, then
 $$\sigma_{p}(M_wEM_u)=\sigma_{jp}(M_wEM_u).$$

\vspace*{0.3cm} {\bf Proof, }
This is a direct consequence of theorems 2.8 and 3.4 of \cite{e2}.\\

\vspace*{0.4cm} {\bf Corollary 3.3.}
If $M_wEM_u$ is quasi-$\ast$-paranormal  and $G=X$ or is an $A$-class operator and $S(E(u))=X$ or is a quasi-$\ast$-$A$-class operator, then $\sigma_{p}(M_wEM_u)=\sigma_{jp}(M_wEM_u)$.\\

\vspace*{0.3cm} {\bf Proof, } This is a direct consequence of Proposition 3.1,  Theorems 2.3 and 2.4.\\

For all $f\in L^2(\Sigma)$ we have\\
\begin{align*}
\langle wE(uf), f\rangle&=\int_{X}wE(uf)\bar{f}d\mu\\
&=\int_{X}ufE(w\bar{f})d\mu\\
&=\int_{X}f\overline{\bar{u}E(\bar{w}f)}d\mu\\
&=\langle f, \bar{u}E(\bar{w}f)\rangle.\\
\end{align*}

This implies that $(M_wEM_u)^{\ast}=M_{\bar{u}}EM_{\bar{w}}$.\\

\vspace*{0.4cm} {\bf Proposition 3.4.} Let $|E(uw)|^2\geq E(|u|^2)E(|w|^2)$ on $G$. Then\\

 (a) Every non-zero isolated point of $ess \
range(E(uw))$ is a simple pole of the resolvent of $M_wEM_u$.\\

 (b) If $\mu$ is a non-zero isolated point of $ess \
range(E(uw))$ and $E_{\mu}$ is the Riesz idempotent of $M_wEM_u$ with respect to $\mu$. Then $E_{\mu}$ is self-adjoint if and only if $N(M_wEM_u-\mu)\subseteq N(M_{\bar{u}}EM_{\bar{w}}-\bar{\mu})$.\\

 \vspace*{0.3cm} {\bf Proof,} By using the results of \cite{sm} and Theorem 2.3 we get the proof.\\

The next corollary is a direct consequence of Theorem  2.6 and the results of \cite{zz}.\\

\vspace*{0.4cm} {\bf Corollary 3.5.} Let $M_wEM_u$ be a bounded operator on $L^2(\Sigma)$. If $\alpha_{\mu}\geq0$ for all $\mu>0$, where
\begin{align*}
\alpha_k=|E(uw)|^{2(n+k)}E(|w|^2)&-(1+n)\mu^n|E(uw)|^{2(k-1)}\chi_{S_0}(E(|w|^2))^2E(|u|^2)\\
&+\mu^{1+n}|E(uw)|^{2(k-1)}\chi_{S_0}E(|w|^2),
\end{align*}
and $S_0=S(E(uw))$,\\
then\\

(a) For every $0\neq\lambda\in \mathbb{C}$ we have
 $$ker (M_wEM_u-\lambda)\subseteq ker(M_{\bar{u}}EM_{\bar{w}}-\bar{\lambda}).$$

(b)
$$\sigma_{jp}(M_wEM_u)\setminus \{0\}=\sigma_p(M_wEM_u)\setminus \{0\},$$
and
$$ \sigma_{ja}(M_wEM_u)\setminus \{0\}=\sigma_a(M_wEM_u)\setminus \{0\}.$$

(c) If $\lambda\neq \mu$, then
$$ker(M_wEM_u-\lambda)\perp ker(M_wEM_u-\mu).$$

(d) For every $0\neq\lambda\in \mathbb{C}$ we have
$$ker(M_wEM_u-\lambda)=ker(M_wEM_u-\lambda)^2,$$
and
$$ ker((M_wEM_u)^{k+1})=ker((M_wEM_u)^{k+2}).$$

(e) The operator $M_wEM_u$ has SVEP.\\

%If $S_1=S(u(E(|w|^2))^{\frac{1}{2}}$, then $L^2(X\backslash
%S_1)\subseteq N(T)$. Also, if
%$S_2=S(\bar(w)(E(|u|^2))^{\frac{1}{2}}$, then $L^2(X\backslash
%S_2)\subseteq N(T^{\ast})$. Thus, if $S_2\neq X$, then $T$ isn't
%dense range. Let $P$ be an orthogonal projection onto
%$$\overline{R(T)}$. Hence $PT=T$, $T^{\ast}P=T^{\ast}$, $PTP=TP$
%and $PT^{\ast}=PT^{\ast}P$.

%The operator $T$ is $M$-paranormal, if for all unit vectors $x$ in
%$H$, $\|Tx\|^2\leq M\|T^2x\|$. $T$ is paranormal if and only if
%$T^{\ast^2}T^2-2kT^{\ast}T+k^2\geq0$, $k\in \mathbb{R}$ and $T$ is
%$M$-paranormal if
%and only if $M^2T^{\ast^2}T^2+2kT^{\ast}T+k^2\geq0$.\\

%If $T$ is paranormal, then
%$$E(|w|^2)|E(u)|^2\leq |E(uw)|(E(|w|^2))^{\frac{1}{2}}|E(u)|.$$

\end{document}